# On the definition, stationary distribution and second order structure of positive semidefinite Ornstein–Uhlenbeck type processes

CHRISTIAN PIGORSCH[1] and ROBERT STELZER[2]

[1]*Institute of Econometrics and Operations Research, Department of Economics, University of Bonn, Adenauerallee 24-42, D-53113 Bonn, Germany.* E-mail: *christian.pigorsch@uni-bonn.de*
[2]*Zentrum Mathematik & TUM Institute for Advanced Study, Technische Universität München, Boltzmannstraße 3, D-85747 Garching, Germany.* E-mail: *rstelzer@ma.tum.de*, url: *http://www-m4.ma.tum.de*

Several important properties of positive semidefinite processes of Ornstein–Uhlenbeck type are analysed. It is shown that linear operators of the form $X \mapsto AX + XA^{\mathrm{T}}$ with $A \in M_d(\mathbb{R})$ are the only ones that can be used in the definition provided one demands a natural non-degeneracy condition. Furthermore, we analyse the absolute continuity properties of the stationary distribution (especially when the driving matrix subordinator is the quadratic variation of a $d$-dimensional Lévy process) and study the question of how to choose the driving matrix subordinator in order to obtain a given stationary distribution. Finally, we present results on the first and second order moment structure of matrix subordinators, which is closely related to the moment structure of positive semidefinite Ornstein–Uhlenbeck type processes. The latter results are important for method of moments based estimation.

## 1. Introduction

Positive semidefinite Ornstein–Uhlenbeck type processes are a special case of Ornstein–Uhlenbeck (OU) type processes (see, e.g., [15, 32]) that take on only values in the positive semidefinite matrices.

They were first introduced in [6] as the solution of the stochastic differential equation

$$\mathrm{d}\Sigma_t = (A\Sigma_{t-} + \Sigma_{t-}A^{\mathrm{T}})\,\mathrm{d}t + \mathrm{d}L_t \qquad (1.1)$$

with $A$ being a $d \times d$-matrix, $L$ a matrix subordinator (see [3]) and the initial value $\Sigma_0$ a positive semidefinite matrix. Matrix subordinators are a generalisation of the one-dimensional concept of Lévy subordinators (see [30]) and are simply Lévy processes that take on only values in the positive semidefinite matrices. As a stochastic process in the positive semidefinite matrices, such an OU type process can be used to describe the random evolution of a covariance matrix over time (that is why we denote it by $\Sigma$). In particular, it can be used to define a multidimensional extension of the popular stochastic volatility model of Barndorff-Nielsen and Shephard [4] (see [25] for details).

In this paper we analyse several important properties of positive semidefinite OU type processes. After briefly reviewing these processes in Section 2, we start in Section 3 by showing that linear operators $\mathbf{A}$ on the symmetric matrices of the form $X \mapsto AX + XA^{\mathrm{T}}$ with a $d \times d$ matrix $A$ are the only ones whose exponential $e^{\mathbf{A}t}$ maps the positive semidefinite matrices onto the positive semidefinite matrices at all times $t \in \mathbb{R}^+$. This implies that these operators are *the* natural class to be used in the definition of positive semidefinite OU type processes. Thereafter we analyse absolute continuity properties of the stationary distribution in Section 4.1. In particular, we show that the stationary distribution is absolutely continuous and, hence, concentrated on the strictly positive definite matrices if one uses the discrete quadratic variation of a $d$-dimensional Lévy processes with absolutely continuous Lévy measure as the driving matrix subordinator (whose Lévy measure is thus concentrated on the rank one matrices). Moreover, we generalise and refine a result of [1] for multivariate OU type processes in Section 4.2. It allows one to compute the characteristic function of a matrix subordinator leading to a given stationary distribution of the OU type process. In connection to this we show also that the stationary solution of (1.1) may be living on the positive semidefinite matrices, although the process $L$ is not a matrix subordinator. This is a very interesting difference from the univariate case, where a stationary OU type process can only be positive at all times if it is driven by a subordinator.

Finally, we aim in Section 5 at characterising the first and second order moment structure (i.e., the mean and covariance matrix) of matrix subordinators in general. Since the first and second order moment structure of OU type processes is a function of the mean and covariance matrix of the driving Lévy process and the drift parameter (i.e., the matrix $A$ above), this characterisation is especially important when one wants to (semi-parametrically) estimate the drift and driving matrix subordinator of a positive semidefinite OU type process or a related stochastic volatility model using a method of moments based estimator. Such estimation techniques are very popular in the univariate case (see, e.g., [4]) and have been successfully used in the multivariate case in [25], for which our results provide the necessary theoretical foundations. We completely characterise in linear algebraic terms the possible first and second moments of diagonal matrix subordinators, i.e., matrix subordinators whose off-diagonal elements are zero. This immediately translates to the first and second moments of the diagonal elements of general matrix subordinators, for which we so far have not been able to derive a complete characterisation of the possible first and second moments. Finally, we show that for the discrete quadratic variation of certain classes of $d$-dimensional Lévy processes, which are of a normal mixture type, the first and second moment can be calculated explicitly in general and that they have a form particularly suitable for further analysis.

## 1.1. Notation

Throughout this paper we write $\mathbb{R}^+$ for the positive real numbers including zero, $\mathbb{R}^{++}$ when zero is excluded and we denote the set of real $n \times n$ matrices by $M_n(\mathbb{R})$, the group of invertible $n \times n$ matrices by $GL_n(\mathbb{R})$, the linear subspace of symmetric matrices by $\mathbb{S}_n$, the (closed) positive semidefinite cone by $\mathbb{S}_n^+$ and the open (in $\mathbb{S}_n$) positive definite cone by $\mathbb{S}_n^{++}$. $I_n$ stands for the $n \times n$ identity matrix and $\sigma(A)$ for the spectrum (the set of all eigenvalues) of a matrix $A \in M_n(\mathbb{R})$. The natural ordering on the symmetric $n \times n$ matrices is denoted by $\leq$, that is, for $A, B \in \mathbb{S}_n$ we have that $A \leq B$, if and only if $B - A \in \mathbb{S}_n^+$. The tensor (Kronecker) product of two matrices $A, B$ is written as $A \otimes B$. vec denotes the well-known vectorisation operator that maps the $n \times n$ matrices to $\mathbb{R}^{n^2}$ by stacking the columns of a matrix below one another. For more information regarding the tensor product and vec operator we refer to [13], Chapter 4. Likewise vech$:\mathbb{S}_d \to \mathbb{R}^{d(d+1)/2}$ denotes the "vector-half" operator that stacks the columns of the lower triangular part of a symmetric matrix below one another. The standard basis matrices (i.e., the matrices which have only zero entries except for a single one in the $i$th row and $j$th column) of $M_d(\mathbb{R})$ are denoted by $E^{(ij)}$ for $i,j = 1, 2, \ldots, d$. Finally, $A^\mathrm{T}$ is the transpose of a matrix $A \in M_n(\mathbb{R})$. For a matrix $A$ we denote by $A_{ij}$ the element in the $i$th row and $j$th column and this notation is extended to processes in a natural way.

Regarding all random variables and processes we assume that they are defined on a given appropriate filtered probability space $(\Omega, \mathcal{F}, P, (\mathcal{F}_t))$ satisfying the usual hypotheses.

Norms of vectors or matrices are denoted by $\|\cdot\|$. If the norm is not specified further, then it is irrelevant which particular norm is used.

Furthermore, we employ an intuitive notation with respect to the (stochastic) integration with matrix-valued integrators referring to any of the standard texts (e.g., [27]) for a comprehensive treatment of the theory of stochastic integration. Let $(L_t)_{t \in \mathbb{R}^+}$ in $M_{n,r}(\mathbb{R})$ be a semimartingale and $(A_t)_{t \in \mathbb{R}^+}$ in $M_{m,n}(\mathbb{R})$, $(B_t)_{t \in \mathbb{R}^+}$ in $M_{r,s}(\mathbb{R})$ be adapted integrable (w.r.t. $L$) processes. Then we denote by $\int_0^t A_s \, \mathrm{d}L_s B_s$ the matrix $C_t$ in $M_{m,s}(\mathbb{R})$ which has $ij$th element

$$C_{ij,t} = \sum_{k=1}^n \sum_{l=1}^r \int_0^t A_{ik,s} B_{lj,s} \, \mathrm{d}L_{kl,s}.$$

Equivalently, such an integral can be understood in the sense of [23], resp. [22], by identifying it with the integral $\int_0^t \mathbf{A}_s \, \mathrm{d}L_s$, with $\mathbf{A}_t$ being for each fixed $t$ the linear operator $M_{n,r}(\mathbb{R}) \to M_{m,s}(\mathbb{R}), X \mapsto A_t X B_t$. Moreover, we always denote by $\int_a^b$ with $a \in \mathbb{R} \cup \{-\infty\}, b \in \mathbb{R}$ the integral over the half-open interval $]a, b]$ for notational convenience. If $b = \infty$ the integral is understood to be over $]a, b[$.

For a set $A$ the indicator function is denoted by $I_A$ and the function $\log^+$ is defined as $\log^+(x) = \max(\log(x), 0)$.

## 2. Positive semidefinite OU type processes

In this section we briefly review the definition of positive semidefinite OU type processes (see [6]) and some of their elementary properties which are relevant for our results later on.

The construction of these processes is based on a special type of matrix-valued Lévy process (see, e.g., [30] for a comprehensive exposition of Lévy processes in general), viz. matrix subordinators, studied in detail in [3]. An $\mathbb{S}_d$-valued Lévy process $L = (L_t)_{t \in \mathbb{R}^+}$ is said to be a *matrix subordinator* if $L_t - L_s \in \mathbb{S}_d^+$ for all $s, t \in \mathbb{R}^+$ with $t > s$. It is $\mathbb{S}_d^+$-increasing and of finite variation.

Let $L$ be a $d \times d$ matrix subordinator, $A \in M_d(\mathbb{R})$ and $\Sigma_0 \in \mathbb{S}_d^+$. Then the solution to the stochastic differential equation of OU type

$$d\Sigma_t = (A\Sigma_{t-} + \Sigma_{t-}A^{\mathrm{T}})\,dt + dL_t$$

is given by

$$\Sigma_t = e^{At}\Sigma_0 e^{A^{\mathrm{T}}t} + \int_0^t e^{A(t-s)}\,dL_s e^{A^{\mathrm{T}}(t-s)}.$$

Note that the drift term simply is the linear operator $\mathbf{A}: \mathbb{S}_d \to \mathbb{S}_d, X \mapsto AX + XA^{\mathrm{T}}$ applied to $\Sigma_{t-}$. Since the increments of $L$ are positive semidefinite and it holds that $e^{\mathbf{A}t}(\mathbb{S}_d^+) = \mathbb{S}_d^+$ for all $t \in \mathbb{R}$, the solution stays in $\mathbb{S}_d^+$ for all times and is therefore called a positive semidefinite OU type process. Provided $E(\log^+ \|L_1\|) < \infty$ and $\sigma(A) \subset (-\infty, 0) + i\mathbb{R}$, the above SDE has a unique stationary solution representable as

$$\Sigma_t = \int_{-\infty}^t e^{A(t-s)}\,dL_s e^{A^{\mathrm{T}}(t-s)}.$$

Next we provide a characterisation of the stationary distribution of the positive semidefinite OU type process. To this end observe that $\mathrm{tr}(X^{\mathrm{T}}Y)$ (with $X, Y \in M_d(\mathbb{R})$ and tr denoting the usual trace functional) defines a scalar product on $M_d(\mathbb{R})$ (respectively, $\mathbb{S}_d$). Moreover, note that the vec operator links the scalar product on $M_d(\mathbb{R})$ ($\mathbb{S}_d$) with the scalar product on $\mathbb{R}^{d^2}$ via $\mathrm{tr}(X^{\mathrm{T}}Y) = \mathrm{vec}(X)^{\mathrm{T}} \mathrm{vec}(Y)$ and that the norm on $M_d(\mathbb{R})$ induced by this scalar product is the Frobenius norm.

This, in particular, implies that the driving matrix subordinator $L$ at time $t \in \mathbb{R}^+$ has characteristic function (cf. [3])

$$\mu_{L_t}(Z) = \exp(t\psi_L(Z)), \qquad (2.1)$$

$$Z \in \mathbb{S}_d, \text{ where } \psi_L(Z) := i\,\mathrm{tr}(\gamma_L Z) + \int_{\mathbb{S}_d^+ \setminus \{0\}} (e^{i\,\mathrm{tr}(XZ)} - 1)\nu_L(dX)$$

with $\gamma_L \in \mathbb{S}_d^+$ and $\nu_L$ being the Lévy measure, which is concentrated on $\mathbb{S}_d^+$ and satisfies

$$\int_{\mathbb{S}_d^+} \min(\|x\|, 1)\nu_L(dx) < \infty.$$

Provided $E(\log^+ \|L_1\|) < \infty$ and $\sigma(A) \subset (-\infty, 0) + i\mathbb{R}$, the stationary distribution of the positive semidefinite Ornstein–Uhlenbeck process $\Sigma$ is again infinitely divisible with characteristic function

$$\hat{\mu}_\Sigma(Z) = \exp\left(\int_0^\infty \psi_L(e^{A^T s} Z e^{As}) \, ds\right)$$
$$= \exp\left(i \operatorname{tr}(\gamma_\Sigma Z) + \int_{\mathbb{S}_d^+ \setminus \{0\}} (e^{i \operatorname{tr}(XZ)} - 1) \nu_\Sigma(dX)\right), \qquad Z \in \mathbb{S}_d, \quad (2.2)$$

where

$$\gamma_\Sigma = -\mathbf{A}^{-1} \gamma_L \quad (2.3)$$

and

$$\nu_\Sigma(E) = \int_0^\infty \int_{\mathbb{S}_d^+ \setminus \{0\}} I_E(e^{As} x e^{A^T s}) \nu_L(dx) \, ds \quad (2.4)$$

for all Borel sets $E$ in $\mathbb{S}_d^+ \setminus \{0\}$ with $I_E$ denoting the indicator function of the set $E$.

Assume that the driving Lévy process is square-integrable. Then the stationary second order moment structure is given by

$$E(\Sigma_t) = \gamma_\Sigma - \mathbf{A}^{-1} \int_{\mathbb{S}_d^+ \setminus \{0\}} y \nu_L(dy) = -\mathbf{A}^{-1} E(L_1), \quad (2.5)$$

$$\operatorname{var}(\operatorname{vec}(\Sigma_t)) = \int_0^\infty e^{(A \otimes I_d + I_d \otimes A)t} \operatorname{var}(\operatorname{vec}(L_1)) e^{(A^T \otimes I_d + I_d \otimes A^T)t} \, dt$$
$$= -\mathcal{A}^{-1} \operatorname{var}(\operatorname{vec}(L_1)), \quad (2.6)$$

$$\operatorname{cov}(\operatorname{vec}(\Sigma_{t+h}), \operatorname{vec}(\Sigma_t)) = e^{(A \otimes I_d + I_d \otimes A)h} \operatorname{var}(\operatorname{vec}(\Sigma_t)), \quad (2.7)$$

where $t \in \mathbb{R}$ and $h \in \mathbb{R}^+$ and $\mathcal{A} : M_{d^2}(\mathbb{R}) \to M_{d^2}(\mathbb{R}), X \mapsto (A \otimes I_d + I_d \otimes A)X + X(A^T \otimes I_d + I_d \otimes A^T)$.

## 3. The possible linear operators

So far linear operators $\mathbf{A} : \mathbb{S}_d \to \mathbb{S}_d$ of the form $X \mapsto AX + XA^T$ with $A \in M_d(\mathbb{R})$ have been used seemingly ad hoc in the definition of positive semidefinite OU type processes. Obviously one can use any linear operator $\mathbf{A}$ which satisfies $e^{\mathbf{A}t}(\mathbb{S}_d^+) \subseteq \mathbb{S}_d^+$ and one would thus like to completely characterise these linear operators. However, this appears unfortunately impossible, because it is a very intricate unsolved problem to characterise all linear operators mapping $\mathbb{S}_d^+$ to $\mathbb{S}_d^+$ but not onto (cf. [16, 24]). In the following we show, however, that linear operators of the form $X \mapsto AX + XA^T$ are the only possible ones if one demands $e^{\mathbf{A}t}(\mathbb{S}_d^+) = \mathbb{S}_d^+$. The latter condition ensures that we have a non-degenerate situation in the following sense: On the one hand, it means that $\Sigma_t$ can take any value in $\mathbb{S}_d^+$ by changing the initial value $\Sigma_0$ as long as the Lévy process does not jump. On the other hand, it obviously implies that the stationary distribution of a positive semidefinite OU type process cannot be concentrated on a proper subset of $\mathbb{S}_d^+$ for a driftless Lévy process unless the Lévy measure of the matrix subordinator is concentrated on a proper subset (combine (2.4) and arguments analogous to [30], Theorem 24.10). Hence, the linear operators of the form $X \mapsto AX + XA^{\mathrm{T}}$ are a natural class to use in the definition of positive semidefinite OU type processes.

Before giving a technical lemma from which our result will follow, it is important to note that the linear operator $\mathbf{A}$ is uniquely characterised by the matrix $A$.

**Proposition 3.1.** *Let $A, C \in M_d(\mathbb{R})$. Then the linear operators $\mathbf{A} : \mathbb{S}_d \to \mathbb{S}_d, X \mapsto AX + XA^{\mathrm{T}}$ and $\mathbf{C} : \mathbb{S}_d \to \mathbb{S}_d, X \mapsto CX + XC^{\mathrm{T}}$ are the same, if and only if $A = C$.*

*Moreover, for any operator $\mathbf{A} : \mathbb{S}_d \to \mathbb{S}_d, X \mapsto AX + XA^{\mathrm{T}}$ the matrix $A$ is already uniquely identified by the values $\{\mathbf{A} E^{(ii)}\}_{i=1,\ldots,d}$.*

**Proof.** It suffices to show the second claim that $\{\mathbf{A} E^{(ii)}\}_{i=1,\ldots,d}$ already uniquely characterises $A$. It is easy to see that

$$(AE^{(ii)} + E^{(ii)} A^{\mathrm{T}})_{kl} = \begin{cases} 2A_{ii}, & \text{for } k = l = i, \\ A_{il}, & \text{for } k = i \text{ and } k \neq l, \\ A_{ik}, & \text{for } l = i \text{ and } k \neq l, \\ 0, & \text{otherwise.} \end{cases}$$

Thus $\mathbf{A} E^{(ii)}$ uniquely characterises the $i$th column of $A$ and, hence, $\{\mathbf{A} E^{(ii)}\}_{i=1,\ldots,d}$ uniquely characterises $A$. $\square$

**Lemma 3.2.** *Let $\mathbf{A} : \mathbb{S}_d \to \mathbb{S}_d$ be a linear operator and assume that there exists an $\varepsilon > 0$ and a function $D : [0, \varepsilon[ \to M_d(\mathbb{R})$ such that $\mathrm{e}^{\mathbf{A} t} Z = D(t) Z D(t)^{\mathrm{T}}$ for all $t \in [0, \varepsilon[$ and $Z \in \mathbb{S}_d$. Then there exists an $\widetilde{\varepsilon} > 0$, a continuously differentiable function $\widetilde{D} : [0, \varepsilon[ \to M_d(\mathbb{R})$ and a unique matrix $A \in M_d(\mathbb{R})$ such that $\mathrm{e}^{\mathbf{A} t} X = \widetilde{D}(t) X \widetilde{D}(t)^{\mathrm{T}}$ for all $t \in [0, \widetilde{\varepsilon}[$ and $X \in \mathbb{S}_d$ and such that $\mathbf{A} X = AX + XA^{\mathrm{T}}$ for all $X \in \mathbb{S}_d$.*

**Proof.** We first show the existence of a matrix $\widetilde{D}(t)$ with the stated properties. Observe that we obtain as a side result that, provided $(\exp(\mathbf{A} t) E^{(11)})_{11} \neq 0$, the operator $\exp(\mathbf{A} t)$ is already fully identified by the values $\exp(\mathbf{A} t) E^{(11)}$ and $\exp(\mathbf{A} t)(E^{(1,j)} + E^{(j,1)})$ with $j = 2, 3, \ldots, d$ and the fact that it can be represented as $X \mapsto D(t) X D(t)^{\mathrm{T}}$.

Note that $D(t) X D(t)^{\mathrm{T}} = (-D(t)) X (-D(t))^{\mathrm{T}}$ for all $X \in \mathbb{S}_d$, so the matrix $D(t)$ can only be unique up to a multiplication by minus one. Elementary calculations give

$$\exp(\mathbf{A} t) E^{(11)} = \begin{pmatrix} D_{11}^2 & D_{11} D_{21} & D_{11} D_{31} & \cdots & D_{11} D_{d1} \\ * & * & * & \cdots & * \end{pmatrix} \quad (3.1)$$

and

$$\exp(\mathbf{A} t)(E^{(1j)} + E^{(j1)})$$

$$= \begin{pmatrix} D_{1j}D_{11} + D_{11}D_{1j} & D_{1j}D_{21} + D_{11}D_{2j} & \cdots & D_{1j}D_{d1} + D_{11}D_{dj} \\ * & * & \cdots & * \end{pmatrix} \quad (3.2)$$

for $j = 2, 3, \ldots, d$. Here $D_{kl}$ denotes $D_{kl}(t)$ for notational convenience and $*$ represents entries which are of no interest in the following. It is easy to see that the above equations uniquely characterise the matrix $D(t)$ up to the sign of $D_{11}(t)$, as long as $(\exp(\mathbf{A}t)E^{(11)})_{11} = D_{11}(t)^2 \neq 0$.

Since $\exp(\mathbf{A} \cdot 0)$ is the identity on $\mathbb{S}_d$ and thus $(\exp(\mathbf{A} \cdot 0)E^{(11)})_{11} = 1$, the continuity of $t \mapsto \exp(\mathbf{A}t)$ ensures that there exists an $\widetilde{\varepsilon} > 0$ with $\widetilde{\varepsilon} \leq \varepsilon$ such that $(\exp(\mathbf{A}t)E^{(11)})_{11} > 0$ for all $t \in [0, \widetilde{\varepsilon}[$. Thus, it follows from (3.1) and (3.2) that $\widetilde{D}(t) \in M_d(\mathbb{R})$ defined by

$$\widetilde{D}_{11}(t) = \sqrt{(\exp(\mathbf{A}t)E^{(11)})_{11}}, \quad (3.3)$$

$$\widetilde{D}_{i1}(t) = \frac{(\exp(\mathbf{A}t)E^{(11)})_{i1}}{\widetilde{D}_{11}(t)} \quad \text{for } i = 2, 3, \ldots, d, \quad (3.4)$$

$$\widetilde{D}_{1j}(t) = \frac{(\exp(\mathbf{A}t)(E^{(1j)} + E^{(j1)}))_{11}}{2\widetilde{D}_{11}(t)} \quad \text{for } j = 2, 3, \ldots, d, \quad (3.5)$$

$$\widetilde{D}_{ij}(t) = \frac{(\exp(\mathbf{A}t)(E^{(1j)} + E^{(j1)}))_{1i}}{\widetilde{D}_{11}(t)} - \frac{(\exp(\mathbf{A}t)(E^{(1j)} + E^{(j1)}))_{11}(\exp(\mathbf{A}t)E^{(11)})_{i1}}{2\widetilde{D}_{11}(t)^3}$$

$$\text{for } i, j = 2, 3, \ldots, d \quad (3.6)$$

is well-defined for all $t \in [0, \widetilde{\varepsilon}[$ and satisfies

$$\exp(\mathbf{A}t)Z = \widetilde{D}(t)Z\widetilde{D}(t)^{\mathrm{T}} \qquad \forall t \in [0, \widetilde{\varepsilon}[ \text{ and } Z \in \mathbb{S}_d.$$

The continuous differentiability of $t \mapsto \exp(\mathbf{A}t)$ and the one of the square root function on $\mathbb{R}^+ \setminus \{0\}$ imply together with (3.3)–(3.6) and the strict positivity of $(\exp(\mathbf{A} \cdot t)E^{(11)})_{11}$ that the map $[0, \widetilde{\varepsilon}[ \to M_d(\mathbb{R}), t \mapsto \widetilde{D}(t)$ is continuously differentiable.

Using the notion of total or Fréchet derivatives (see [29] or [9], Section X.4, for a review in connection with matrix analysis) and denoting the linear operators on $\mathbb{S}_d$ by $L(\mathbb{S}_d)$, it follows immediately from

$$(X + H)Z(X + H)^{\mathrm{T}} = XZX^{\mathrm{T}} + XZH^{\mathrm{T}} + HZX^{\mathrm{T}} + HZH^{\mathrm{T}} \qquad \forall Z \in \mathbb{S}_d$$

for $X, H \in M_d(\mathbb{R})$ that the map $f : M_d(\mathbb{R}) \to L(\mathbb{S}_d), X \mapsto f(X)$ with $f(X)Z = XZX^{\mathrm{T}}$ for all $Z \in \mathbb{S}_d$ is continuously differentiable and the derivative $\mathcal{D}f(X)$ is given by the linear map $M_d(\mathbb{R}) \mapsto L(\mathbb{S}_d), H \mapsto \mathcal{D}f(X)(H)$ with $\mathcal{D}f(X)(H)Z = XZH^{\mathrm{T}} + HZX^{\mathrm{T}}$ for $Z \in \mathbb{S}_d$. Thus

$$\left(\frac{\mathrm{d}}{\mathrm{d}t}\exp(\mathbf{A}t)\right)Z = \left(\frac{\mathrm{d}}{\mathrm{d}t}f(\widetilde{D}(t))\right)Z = \mathcal{D}f(\widetilde{D}(t))\left(\frac{\mathrm{d}}{\mathrm{d}t}\widetilde{D}(t)\right)Z$$

$$= \widetilde{D}(t)Z\left(\frac{\mathrm{d}}{\mathrm{d}t}\widetilde{D}(t)\right)^{\mathrm{T}} + \left(\frac{\mathrm{d}}{\mathrm{d}t}\widetilde{D}(t)\right)Z\widetilde{D}(t)^{\mathrm{T}}$$

for all $t \in [0, \tilde{\varepsilon}[$ and $Z \in \mathbb{S}_d$.

Since $\frac{\mathrm{d}}{\mathrm{d}t}\exp(\mathbf{A}t) = \exp(\mathbf{A}t)\mathbf{A}$, it follows that $\mathbf{A} = \exp(-\mathbf{A}t)\frac{\mathrm{d}}{\mathrm{d}t}\exp(\mathbf{A}t)$. Moreover, it is easy to see that $\widetilde{D}(t) \in GL_d(\mathbb{R})$ and $\exp(\mathbf{A}t)^{-1}Z = \exp(-\mathbf{A}t)Z = \widetilde{D}(t)^{-1}Z\widetilde{D}(t)^{-\mathrm{T}}$ for $Z \in \mathbb{S}_d$ and hence

$$\mathbf{A}Z = \widetilde{D}(t)^{-1}\left(\widetilde{D}(t)Z\left(\frac{\mathrm{d}}{\mathrm{d}t}\widetilde{D}(t)\right)^{\mathrm{T}} + \left(\frac{\mathrm{d}}{\mathrm{d}t}\widetilde{D}(t)\right)Z\widetilde{D}(t)^{\mathrm{T}}\right)\widetilde{D}(t)^{-\mathrm{T}} \quad (3.7)$$

$$= Z\left(\widetilde{D}(t)^{-1}\frac{\mathrm{d}}{\mathrm{d}t}\widetilde{D}(t)\right)^{\mathrm{T}} + \left(\widetilde{D}(t)^{-1}\frac{\mathrm{d}}{\mathrm{d}t}\widetilde{D}(t)\right)Z \quad (3.8)$$

for $Z \in \mathbb{S}_d$ and all $t \in [0, \tilde{\varepsilon}[$. As by construction $\widetilde{D}(0) = I_d$, setting $A = \frac{\mathrm{d}}{\mathrm{d}t}\widetilde{D}(t)|_{t=0}$ concludes the proof now by noting that Proposition 3.1 ensures the uniqueness of $A \in M_d(\mathbb{R})$. □

Using this result, we can now completely characterise the linear operators whose exponential group maps the positive semidefinite matrices onto themselves, as desired.

**Theorem 3.3.** *Let $d \in \mathbb{N}$. Then the following holds:*

(i) *A linear mapping $\mathbf{A}: \mathbb{S}_d \to \mathbb{S}_d$ satisfies $\mathrm{e}^{\mathbf{A}t}(\mathbb{S}_d^+) = \mathbb{S}_d^+$ for all $t \in \mathbb{R}^+$, if and only if there exists a matrix $A \in M_d(\mathbb{R})$ such that $\mathbf{A}X = AX + XA^{\mathrm{T}}$ for all $X \in \mathbb{S}_d$.*

(ii) *A linear mapping $\mathbf{A}: \mathbb{S}_d \to \mathbb{S}_d$ satisfies $\mathrm{e}^{\mathbf{A}t}(\mathbb{S}_d^+) = \mathbb{S}_d^+$ for all $t \in \mathbb{R}$, if there exists an $\varepsilon > 0$ such that $\mathrm{e}^{\mathbf{A}t}(\mathbb{S}_d^+) = \mathbb{S}_d^+$ for all $t \in [0, \varepsilon[$.*

**Proof.** (i) Note that all linear maps that map $\mathbb{S}_d^+$ onto itself are of the form $X \mapsto CXC^{\mathrm{T}}$ for some $C \in GL_d(\mathbb{R})$ (cf. [24], Chapter 3, or the original article [33]). The "if" part follows by combining this with the fact that $\mathrm{e}^{\mathbf{A}t}X = \mathrm{e}^{At}X\mathrm{e}^{A^{\mathrm{T}}t}$ for all $X \in \mathbb{S}_d$ and $t \in \mathbb{R}^+$ (see, e.g., [13], pp. 255 and 440). Regarding the "only if" part, this shows that there is a function $D(t): \mathbb{R} \to M_d(\mathbb{R})$ such that $\mathrm{e}^{\mathbf{A}t}X = D(t)XD(t)^{\mathrm{T}}$ for all $t \in \mathbb{R}$ and $X \in \mathbb{S}_d$. Lemma 3.2 immediately concludes now.

(ii) Follows immediately from Lemma 3.2 as well. □

$\mathbb{S}_d^+$ can be replaced by $\mathbb{S}_d^{++}$ in the above result. Yet we cannot extend the result to the case $\mathrm{e}^{\mathbf{A}t}(\mathbb{S}_d^+) \subset \mathbb{S}_d^+$, as there are linear operators $\mathbf{C}$ such that $\mathbf{C}(\mathbb{S}_d^+) \subset \mathbb{S}_d^+$ which are not representable by $X \mapsto CXC^{\mathrm{T}}$ for some $C \in M_d(\mathbb{R})$ (cf. [11]).

## 4. The stationary distribution

In this section we first use the theory of operator self-decomposability to analyse the absolute continuity properties of the stationary distribution and then show how to find a driving Lévy process for a given stationary distribution. Finally, we obtain that the stationary distribution – somewhat surprisingly – can be concentrated on $\mathbb{S}_d^+$, although $L$ is not a matrix subordinator.

### 4.1. Operator self-decomposability and absolute continuity

The stationary distributions of multivariate Ornstein–Uhlenbeck type processes are operator self-decomposable (see [32]). The operator self-decomposable distributions (cf. [15] for a comprehensive treatment) form an important subclass of the infinitely divisible distributions and the definition adapted to our matrix case reads as follows:

**Definition 4.1.** *Let $Q:\mathbb{S}_d \to \mathbb{S}_d$ be a linear operator. A probability distribution $\mu$ on $\mathbb{S}_d$ is called* operator self-decomposable with respect to the operator $Q$ *if there exists a probability distribution $\nu_t$ on $\mathbb{S}_d$ such that $\mu = (\mathrm{e}^{Qt}\mu) \ast \nu_t$ for all $t \in \mathbb{R}^+$, where $\ast$ denotes the usual convolution of probability measures.*

**Remark 4.2.** (i) $\nu_t$ is for all $t \in \mathbb{R}^+$ infinitely divisible ([15], Theorem 3.3.5).

(ii) Rephrasing the definition on the level of random variables, we call an $\mathbb{S}_d$-valued random variable $X$ operator self-decomposable with respect to a linear operator $Q:\mathbb{S}_d \to \mathbb{S}_d$ if there exists an $\mathbb{S}_d^+$-valued random variable $Y_t$ independent of $X$ such that $X \stackrel{\text{law}}{=} \mathrm{e}^{Qt}X + Y_t$ for all $t \in \mathbb{R}^+$.

For our Ornstein–Uhlenbeck type processes in the positive semidefinite matrices we immediately have the following result:

**Proposition 4.3.** *The stationary distribution of the positive semidefinite OU type process $\Sigma$ is operator self-decomposable with respect to the operator $\mathbf{A}:\mathbb{S}_d \to \mathbb{S}_d, X \mapsto AX + XA^{\mathrm{T}}$.*

This observation allows us to use the well-established theory on operator self-decomposable distributions to study further properties of our processes. We can now show that the stationary distribution always has a density (with respect to the Lebesgue measure) whenever the Lévy measure of the driving Lévy process has a non-degenerate (with respect to $A$) support. Note that by the Lebesgue measure $\lambda_{\mathbb{S}_d}$ on $\mathbb{S}_d$ we mean the image of the Lebesgue measure $\lambda_{\mathbb{R}^{d(d+1)/2}}$ on $\mathbb{R}^{d(d+1)/2}$ under $\mathrm{vech}:\mathbb{S}_d \to \mathbb{R}^{d(d+1)/2}$, that is, $\lambda_{\mathbb{S}_d}(B) = \lambda_{\mathbb{R}^{d(d+1)/2}}(\mathrm{vech}(B))$ for all Borel sets $B \subset \mathbb{S}_d$.

**Theorem 4.4.** *Assume that the support of the Lévy measure $\nu_L \neq 0$ of the driving matrix subordinator $L$ is contained in a linear subspace $V \subseteq \mathbb{S}_d$, that the smallest $\mathbf{A}$-invariant linear subspace of $\mathbb{S}_d$ containing $V$ is $\mathbb{S}_d$ itself and that $\nu_L(W) = 0$ for every proper linear subspace $W \subset V$.*

*Then the stationary distribution of the Ornstein–Uhlenbeck type process $\Sigma$ is absolutely continuous with respect to the Lebesgue measure on $\mathbb{S}_d$. Moreover, the stationary distribution $P_\Sigma$ of $\Sigma$ is concentrated on $\mathbb{S}_d^{++}$, that is $P_\Sigma(\mathbb{S}_d^{++}) = 1$.*

**Proof.** This follows by a straightforward adaptation of [15], Proposition 3.8.6, or [37] to the $\mathbb{S}_d$-case and the fact that $\lambda_{\mathbb{S}_d}(\mathbb{S}_d^+ \setminus \mathbb{S}_d^{++}) = 0$, which is shown in the next lemma. $\square$

Observe that this result holds for OU type processes in $\mathbb{S}_d$ in general, not only for positive semidefinite ones, and that the density can be characterised as the unique weak solution to a certain integro-differential equation (see [15], Theorem 3.8.10).

**Lemma 4.5.** *The Lebesgue measure on $\mathbb{S}_d$ satisfies $\lambda_{\mathbb{S}_d}(\mathbb{S}_d^+ \setminus \mathbb{S}_d^{++}) = 0$.*

**Proof.** It is clear that $\mathbb{S}_d^+ \setminus \mathbb{S}_d^{++} \subseteq \{X \in \mathbb{S}_d : \det(X) = 0\}$. The determinant is a polynomial in the entries of $X \in \mathbb{S}_d$. Hence, there exists a $d(d+1)/2$-variate polynomial $P : \mathbb{R}^{d(d+1)/2} \to \mathbb{R}$ such that $\det(X) = P(\mathrm{vech}(X))$ and $P$ is not identically zero. As a polynomial $P$ is a real analytic function and the set of zeros of real analytical functions not identically zero has Lebesgue measure zero ([15], Lemma 3.8.4), we have

$$\lambda_{\mathbb{S}_d}(\{X \in \mathbb{S}_d : \det(X) = 0\}) = \lambda_{\mathbb{R}^{d(d+1)/2}}(\mathrm{vech}(\{X \in \mathbb{S}_d : \det(X) = 0\}))$$
$$= \lambda_{\mathbb{R}^{d(d+1)/2}}(\{x \in \mathbb{R}^{d(d+1)/2} : P(x) = 0\}) = 0. \quad \square$$

The following is an immediate consequence of Theorem 4.4 which will be sufficient in most cases encountered in practice.

**Corollary 4.6.** *Assume that for every proper linear subspace $V$ of $\mathbb{S}_d$ the Lévy measure $\nu_L$ of the driving matrix subordinator satisfies $\nu_L(V) = 0$ and that $\nu_L \neq 0$.*

*Then the stationary distribution of the Ornstein–Uhlenbeck type process $\Sigma$ is absolutely continuous with respect to the Lebesgue measure. Moreover, the stationary distribution $P_\Sigma$ of $\Sigma$ is concentrated on $\mathbb{S}_d^{++}$, that is, $P_\Sigma(\mathbb{S}_d^{++}) = 1$.*

A particular class of $d \times d$ matrix subordinators is given by the discrete part of the quadratic variation of a $d$-dimensional Lévy process. Let $\tilde{L}$ thus be a Lévy process in $\mathbb{R}^d$. Then the discrete part of the quadratic variation is given by

$$[\tilde{L}, \tilde{L}]_t^{\mathfrak{d}} := \sum_{0 \leq s \leq t} \Delta \tilde{L}_s (\Delta \tilde{L}_s)^{\mathrm{T}}$$

and it is easy to see that $([\tilde{L}, \tilde{L}]_t^{\mathfrak{d}})_{t \in \mathbb{R}^+}$ is a matrix subordinator in $\mathbb{S}_d$ with the Lévy measure given by $\nu_{[\tilde{L}, \tilde{L}]^{\mathfrak{d}}}(B) = \int_{\mathbb{R}^d} I_B(xx^{\mathrm{T}}) \nu_{\tilde{L}}(\mathrm{d}x) = \nu_{\tilde{L}}(\{x : xx^{\mathrm{T}} \in B\})$ where $\nu_{\tilde{L}}$ denotes the Lévy measure of $\tilde{L}$. Note that $\nu_{[\tilde{L}, \tilde{L}]^{\mathfrak{d}}}$ is obviously concentrated on the positive semidefinite rank one matrices and has thus a rather degenerate support. Yet, under rather mild regularity conditions, the stationary distribution of a positive semidefinite OU type process $\Sigma$ driven by $[\tilde{L}, \tilde{L}]^{\mathfrak{d}}$ is absolutely continuous.

**Proposition 4.7.** *Let $\tilde{L}$ be a Lévy process in $\mathbb{R}^d$ with Lévy measure $\nu_{\tilde{L}} \neq 0$ and assume that $\nu_{\tilde{L}}(\{x : x^{\mathrm{T}} Z x = 0\}) = 0$ for all $Z \in \mathbb{S}_d \setminus \{0\}$.*

*Then the stationary distribution of the Ornstein–Uhlenbeck type process $\Sigma$ driven by $[\tilde{L}, \tilde{L}]^{\diamond}$ is absolutely continuous with respect to the Lebesgue measure. Moreover, the stationary distribution $P_\Sigma$ of $\Sigma$ is concentrated on $\mathbb{S}_d^{++}$, that is, $P_\Sigma(\mathbb{S}_d^{++}) = 1$.*

**Proof.** Recall that $\mathbb{S}_d$ with the scalar product $X, Y \mapsto \operatorname{tr}(X^{\mathrm{T}} Y)$ is a Hilbert space and denote for a set $H \subset \mathbb{S}_d$ its orthogonal complement (with respect to this scalar product) by $H^\perp$.

Let $H$ now be a proper linear subspace of $\mathbb{S}_d$ and $Z \in H^\perp \setminus \{0\}$. Then $H \subseteq \{Z\}^\perp$. Hence, $\nu_{[\tilde{L}, \tilde{L}]^{\diamond}}(H) \leq \nu_{[\tilde{L}, \tilde{L}]^{\diamond}}(\{Z\}^\perp) = \nu_{\tilde{L}}(\{x \in \mathbb{R}^d : xx^{\mathrm{T}} \in \{Z\}^\perp\}) = \nu_{\tilde{L}}(\{x \in \mathbb{R}^d : \operatorname{tr}(xx^{\mathrm{T}} Z) = x^{\mathrm{T}} Z x = 0\}) = 0$ by assumption.

Since $\nu_{[\tilde{L}, \tilde{L}]^{\diamond}} \neq 0$ is obviously true, the proposition is implied by Corollary 4.6. $\square$

The following corollary is usually sufficient for applications:

**Corollary 4.8.** *Let $\tilde{L}$ be a Lévy process in $\mathbb{R}^d$ with Lévy measure $\nu_{\tilde{L}} \neq 0$ and assume that $\nu_{\tilde{L}}$ is absolutely continuous with respect to the Lebesgue measure on $\mathbb{R}^d$.*

*Then the stationary distribution of the Ornstein–Uhlenbeck type process $\Sigma$ driven by $[\tilde{L}, \tilde{L}]^{\diamond}$ is absolutely continuous with respect to the Lebesgue measure. Moreover, the stationary distribution $P_\Sigma$ of $\Sigma$ is concentrated on $\mathbb{S}_d^{++}$, that is, $P_\Sigma(\mathbb{S}_d^{++}) = 1$.*

**Proof.** Let $Z \in \mathbb{S}_d \setminus \{0\}$. Then the function $\mathbb{R}^d \to \mathbb{R}, x = (x_1, x_2, \ldots, x_d)^{\mathrm{T}} \mapsto x^{\mathrm{T}} Z x$ is a quadratic polynomial (in $d$ variables) and thus a real analytic function not identically zero. The zeros of a real analytic function not identically zero form a set of Lebesgue measure zero (see, e.g., [15], Lemma 3.8.4). Hence, $\nu_{\tilde{L}}(\{x : x^{\mathrm{T}} Z x = 0\}) = 0$ for all $Z \in \mathbb{S}_d \setminus \{0\}$ and thus Proposition 4.7 concludes the proof. $\square$

Furthermore, note that the result concerning the existence of an infinitely often differentiable density with all derivatives bounded given in [19], Proposition 3.31, is also immediately applicable to our positive semidefinite OU type processes when using the vech operator to transfer them to $\mathbb{R}^{d(d+1)/2}$-valued OU type processes.

## 4.2. Finding a driving Lévy process

It is also possible to construct positive semidefinite OU type processes with a prescribed marginal stationary distribution. The result is a refined version of [1], Lemma 5.1, adapted to the matrix case focusing especially on matrix subordinators. "Differentiability" in the following means total or Fréchet differentiability and "derivative" refers to the total derivative. Moreover, observe in the following that a probability distribution $\mu$ concentrated on $\mathbb{S}_d^+$ is infinitely divisible if and only if its characteristic function $\hat{\mu}$ is of the form

$$\hat{\mu}(Z) = \exp\biggl(\mathrm{i}\operatorname{tr}(\gamma_\mu Z) + \int_{\mathbb{S}_d} (\mathrm{e}^{\mathrm{i}\operatorname{tr}(XZ)} - 1)\nu_\mu(\mathrm{d}X)\biggr), \qquad Z \in \mathbb{S}_d, \qquad (4.1)$$

where $\gamma_\mu \in \mathbb{S}_d^+$ and $\nu_\mu$ is a Lévy measure on $\mathbb{S}_d$ satisfying

$$\nu_\mu(\mathbb{S}_d \setminus \mathbb{S}_d^+) = 0 \quad \text{and} \quad \int_{\|X\| \leq 1} \|X\| \nu_\mu(\mathrm{d}X) < \infty.$$

This follows from [34], page 156, combined with arguments analogous to the ones in [30], Remark 24.9.

**Theorem 4.9.** *Let $\mathbf{A} : \mathbb{S}_d \to \mathbb{S}_d$ be a linear operator such that there is an $A \in M_d(\mathbb{R})$ with $\sigma(A) \subset (-\infty, 0) + \mathrm{i}\mathbb{R}$ satisfying $\mathbf{A}X = AX + XA^{\mathrm{T}}$ for all $X \in \mathbb{S}_d$. Furthermore, let $\mu$ be an operator self-decomposable (with respect to $\mathbf{A}$) distribution on $\mathbb{S}_d^+$ with associated $\gamma_\mu$ and $\nu_\mu$ such that (4.1) holds and let $\psi(Z) = \log(\int_{\mathbb{S}_d^+} \mathrm{e}^{\mathrm{i}\,\mathrm{tr}(xZ)} \mu(\mathrm{d}x)), Z \in \mathbb{S}_d$, be its cumulant transform (logarithm of the characteristic function $\hat{\mu}$). Assume that $-A\gamma_\mu - \gamma_\mu A^{\mathrm{T}} \in \mathbb{S}_d^+$, that $\psi(Z)$ is differentiable for all $Z \neq 0$ with derivative $D\psi(Z)$ and that*

$$\psi_L : Z \mapsto \begin{cases} -D\psi(Z)(A^{\mathrm{T}}Z + ZA), & \text{for } Z \in \mathbb{S}_d \setminus \{0\}, \\ 0, & \text{for } Z = 0, \end{cases}$$

*is continuous at zero. Then $\hat{\mu}_L : Z \mapsto \exp(\psi_L(Z))$ is the characteristic function of an infinitely divisible distribution $\mu_L$ on $\mathbb{S}_d^+$.*

*Let $L$ be a matrix subordinator with characteristic function $\hat{\mu}_L$ at time one. Then the positive semidefinite Ornstein–Uhlenbeck type process $\mathrm{d}\Sigma_t = \mathbf{A}\Sigma_t\,\mathrm{d}t + \mathrm{d}L_t$ driven by $L$ has stationary distribution $\mu$.*

**Proof.** As the characteristic functions of infinitely divisible distributions have no zeros ([30], Lemma 7.5), $\psi$ is well defined. The definition of operator self-decomposability implies together with Remark 4.2(i) that $\exp(\psi(Z))/\exp(\psi(\mathrm{e}^{A^{\mathrm{T}}t} Z \mathrm{e}^{At}))$ is for all $t \in \mathbb{R}^+$ the characteristic function of an infinitely divisible distribution. Hence the function

$$Z \mapsto (\exp(\psi(Z))/\exp(\psi(\mathrm{e}^{A^{\mathrm{T}}t} Z \mathrm{e}^{At})))^{1/t}$$

is for all $t \in \mathbb{R}^+$ the characteristic function of an infinitely divisible distribution. Since, obviously,

$$\frac{\psi(Z) - \psi(\mathrm{e}^{A^{\mathrm{T}}t} Z \mathrm{e}^{At})}{t} \to -D\psi(Z)(A^{\mathrm{T}}Z + ZA) \qquad \text{as } t \to 0$$

for all $Z \in \mathbb{S}_d$, we have that $(\exp(\psi(Z))/\exp(\psi(\mathrm{e}^{A^{\mathrm{T}}t} Z \mathrm{e}^{At})))^{1/t} \to \exp(\psi_L(Z))$ point-wise for all $Z \in \mathbb{S}_d$ as $t \to 0$. The continuity of $\psi_L$ at zero shows that $\hat{\mu}_L$ is a characteristic function. $\hat{\mu}_L$ belongs to an infinitely divisible distribution on $\mathbb{S}_d$, because infinite divisibility is preserved under weak convergence ([30], Lemma 7.8).

Let $(L_t)_{t \in \mathbb{R}}$ now be a Lévy process with characteristic function $\hat{\mu}_L$ at time one and assume that $\mathrm{d}\Sigma_t = \mathbf{A}\Sigma_t\,\mathrm{d}t + \mathrm{d}L_t$ is started at time zero with $\Sigma_0$ being independent of $(L_t)_{t \in \mathbb{R}}$ and having distribution $\mu$. Thus $\Sigma_t = \mathrm{e}^{At}\Sigma_0 \mathrm{e}^{A^{\mathrm{T}}t} + \int_0^t \mathrm{e}^{A(t-s)}\,\mathrm{d}L_s \mathrm{e}^{A^{\mathrm{T}}(t-s)}$ and we have from [28] (see also [31], [15], Lemma 3.6.4, or the review in [19], Section 2.2) that

$$E(\exp(i\operatorname{tr}(\Sigma_t Z)))$$
$$= \hat{\mu}(e^{A^T t} Z e^{At}) \exp\left(\int_0^t \psi_L(e^{A^T s} Z e^{As})\,ds\right)$$
$$= \hat{\mu}(e^{A^T t} Z e^{At}) \exp\left(\int_0^t -D\psi(e^{A^T s} Z e^{As})(A^T e^{A^T s} Z e^{As} + e^{A^T s} Z e^{As} A)\,ds\right)$$
$$= \hat{\mu}(e^{A^T t} Z e^{At}) \exp(\psi(Z) - \psi(e^{A^T t} Z e^{At})) = \hat{\mu}(Z).$$

Hence, $\Sigma$ is stationary with stationary distribution $\mu$. By [32], Theorem 4.2, it follows now that $E(\log^+ \|L_t\|) < \infty$. As the stationary distribution $\mu$ is concentrated on $\mathbb{S}_d^+$, $\mu_L$ has to be an infinitely divisible distribution on $\mathbb{S}_d^+$ and hence the Lévy process $L$ is a matrix subordinator. This can be seen in the following three steps:

(i) Assume that the Lévy measure $\nu_L$ was not concentrated on $\mathbb{S}_d^+$, that is, $\nu_L(\mathbb{S}_d \setminus \mathbb{S}_d^+) > 0$. Then, using the results of Section 3, which imply for all $s \in \mathbb{R}^+$ that $e^{As} x e^{A^T s} \in \mathbb{S}_d^+$ if and only if $x \in \mathbb{S}_d^+$, and (2.4) (which also holds in this general case), we have that

$$\int_{\mathbb{S}_d} I_{\mathbb{S}_d \setminus \mathbb{S}_d^+}(e^{As} x e^{A^T s}) \nu_L(dx) > 0 \qquad \text{for all } s \in \mathbb{R}^+$$

and hence $\nu_\mu(\mathbb{S}_d \setminus \mathbb{S}_d^+) = \int_0^\infty \int_{\mathbb{S}_d} I_{\mathbb{S}_d \setminus \mathbb{S}_d^+}(e^{As} x e^{A^T s}) \nu_L(dx)\,ds > 0$, a contradiction.

(ii) We have for any $t^* \in \mathbb{R}^{++}$ by (2.4) and (4.1) that

$$\int_0^{t^*} \int_{\mathbb{S}_d^+} (\|e^{As} x e^{A^T s}\| \wedge 1) \nu_L(dx) \leq \int_{\mathbb{S}_d} (\|x\| \wedge 1) \nu_\Sigma(dx) < \infty.$$

Choosing $t^*$ such that $\|e^{As} x e^{A^T s}\| \geq \|x\|/2$ for all $s \in [0, t^*]$ and $x \in \mathbb{S}_d$ implies that

$$\int_0^{t^*} \int_{\mathbb{S}_d^+} (\|e^{As} x e^{A^T s}\| \wedge 1) \nu_L(dx) \geq \int_0^{t^*} \int_{\mathbb{S}_d^+} ((\|x\|/2) \wedge 1) \nu_L(dx)$$
$$= 2t^* \int_{\mathbb{S}_d^+, \|x\| \leq 2} \|x\| \nu_L(dx).$$

Hence, $\int_{\mathbb{S}_d^+, \|x\| \leq 1} \|x\| \nu_L(dx) < \infty$ and therefore $L$ is of finite variation.

(iii) Combining (i) and (ii) implies that $\hat{\mu}_L$ is of the form (4.1) with $\gamma_L \in \mathbb{S}_d$. Yet, (2.3) implies that $\gamma_L = -A\gamma_\mu - \gamma_\mu A^T$ and thus $\gamma_L \in \mathbb{S}_d^+$ by assumption. □

**Remark 4.10.** (a) $D\psi(Z)$ is for each $Z \in \mathbb{S}_d$ a linear operator from $\mathbb{S}_d$ to $\mathbb{R}$. Since $\mathbb{S}_d$ is a Hilbert space, the Riesz–Frechet theorem ensures that there exists a function $\hat{D}_\psi : \mathbb{S}_d \to \mathbb{S}_d$ with $D\psi(Z)X = \text{tr}(\hat{D}_\psi(Z)X)$ for all $X \in \mathbb{S}_d$. This gives an alternative representation of $\psi_L$ above.

(b) The condition $-A\gamma_\mu - \gamma_\mu A^\mathrm{T} \in \mathbb{S}_d^+$ may appear superfluous at a first sight, since it obviously is always satisfied for $d = 1$. However, in general dimensions it need not hold. For example, for

$$A = \begin{pmatrix} -1/10 & -1/3 \\ -1/3 & -2 \end{pmatrix} \quad \text{and} \quad \gamma_\mu = \begin{pmatrix} 2 & -2/3 \\ -2/3 & 2 \end{pmatrix}$$

we have that

$$\sigma(A) = \{-\tfrac{21}{20} \pm \tfrac{1}{60}\sqrt{3649}\} \approx \{-0.043, -2.06\},$$
$$\sigma(\gamma_\mu) = \{8/3, 4/3\}, \quad \text{but } \sigma(-A\gamma_\mu - \gamma_\mu A^\mathrm{T}) = \{\tfrac{169}{45} \pm \tfrac{1}{3}\sqrt{130}\} \approx \{7.56, -0.045\}.$$

If the condition is not satisfied, all assertions of the theorem remain valid except that then $\mu_L$ is not concentrated on $\mathbb{S}_d^+$ and $L$ is a finite variation Lévy process in $\mathbb{S}_d$ with all jumps in $\mathbb{S}_d^+$ but a drift in $\mathbb{S}_d \setminus \mathbb{S}_d^+$.

(c) Choosing $A$ and $\gamma_\mu$ as in (b), setting $\gamma_L = -A\gamma_\mu - \gamma_\mu A^\mathrm{T}$ and choosing $\nu_L$ as some measure on $\mathbb{S}_d^+$ satisfying $\nu_L(\{0\}) = 0$ and $\int_{\mathbb{S}_d^+}(\|x\| \wedge 1)\nu_L(\mathrm{d}x) < \infty$, immediately gives an OU type process in $\mathbb{S}_d$ with stationary distribution concentrated on $\mathbb{S}_d^+$, which is, however, not driven by a matrix subordinator.

## 5. The covariance structure of matrix subordinators

A very important topic is the characterisation of the second order structure of matrix subordinators, that is, which positive semidefinite $d^2 \times d^2$ matrices appear as covariance matrices of $d \times d$ matrix subordinators. Until now we have only been able to find a useful characterisation of the covariance matrices of two subclasses of matrix subordinators, viz. diagonal and special quadratic variation matrix subordinators.

The results obtained below provide insight into the specification and estimation of concrete models in applications. For example, [25] considers a stochastic volatility model and estimates of the mean and covariance matrix of the driving diagonal matrix subordinator of the positive semidefinite OU type stochastic volatility process. For such an estimation procedure we completely characterise the range of possible values for the estimators below.

### 5.1. Diagonal matrix subordinators

Let $L_1, L_2, \ldots, L_d$ be $d$ univariate Lévy subordinators forming together a $d$-dimensional Lévy process. Then the process

$$L = \begin{pmatrix} L_1 & 0 & \cdots & 0 \\ 0 & L_2 & \ddots & \vdots \\ \vdots & \ddots & \ddots & 0 \\ 0 & \cdots & 0 & L_d \end{pmatrix}$$

is a matrix subordinator with non-zero elements only on the diagonal. Such matrix subordinators are referred to as *diagonal matrix subordinators* in the following. Obviously, all information about $L$, in particular the second order structure, is already contained in the diagonal $((L_{1,t}, L_{2,t}, \ldots, L_{d,t})^{\mathrm{T}})_{t \in \mathbb{R}^+}$. It should be noted that $((L_{1,t}, L_{2,t}, \ldots, L_{d,t})^{\mathrm{T}})_{t \in \mathbb{R}^+}$ is increasing in the usual order of the cone $(\mathbb{R}^+)^d$ and so the $d$-dimensional Lévy process $((L_{1,t}, L_{2,t}, \ldots, L_{d,t})^{\mathrm{T}})_{t \in \mathbb{R}^+}$ is a generalisation of the notion of a subordinator to the space $\mathbb{R}^d$. In the following we call a $d$-dimensional Lévy processes with all components being univariate subordinators a *multivariate subordinator*.

For the analysis of the second order structure of multivariate subordinators (and thus of diagonal matrix subordinators) we need the following notions from linear algebra, referring to [7, 8, 36] and the references therein for further details.

**Definition 5.1.** *A matrix $A \in M_d(\mathbb{R})$ is called* completely positive *if there exist a $k \in \mathbb{N}$ and a $B \in M_{d,k}(\mathbb{R})$ with all entries being non-negative such that $A = BB^{\mathrm{T}}$.*

*A matrix $A \in M_d(\mathbb{R})$ is called* doubly non-negative *if $A$ is positive semidefinite and all entries are non-negative.*

Obviously any completely positive matrix is doubly non-negative.

**Proposition 5.2.** (a) *If $L$ is a $d$-dimensional multivariate subordinator with finite second moments, $\mathrm{var}(L_1)$ is completely positive.*

(b) *If $C \in M_d(\mathbb{R})$ is completely positive, there exists a $d$-dimensional multivariate subordinator $L$ such that $C = \mathrm{var}(L_1)$.*

(c) *For $C \in M_d(\mathbb{R})$ there exists a $d$-dimensional multivariate subordinator $L$ with $\mathrm{var}(L_1) = C$ if and only if $C$ is completely positive.*

(d) *Provided $d \leq 4$, there exists a $d$-dimensional multivariate subordinator $L$ with $\mathrm{var}(L_1) = C$ if and only if $C \in M_d(\mathbb{R})$ is doubly non-negative.*

**Proof.** Regarding part (a) observe that $\mathrm{var}(L_1) = \int_{\mathbb{R}^d} xx^{\mathrm{T}} \nu_L(\mathrm{d}x)$ (cf. [30], page 163) and that the Lévy measure $\nu_L$ of $L$ is necessarily concentrated on $(\mathbb{R}^+)^d$. Thus approximating $\int_{(\mathbb{R}^+)^d} xx^{\mathrm{T}} \nu_L(\mathrm{d}x)$ by a sequence of integrals of simple functions $(f_n)_{n \in \mathbb{N}}$ of the form $f_n(x) = \sum_{i=1}^{k} x_{i,n} x_{i,n}^{\mathrm{T}} I_{A_{i,n}}(x)$ with appropriate $x_{i,n} \in (\mathbb{R}^+)^d$ and $A_{i,n} \subset (\mathbb{R}^+)^d$ converging to $x \mapsto xx^{\mathrm{T}}$ gives that $\mathrm{var}(L_1)$ is the limit of a sequence of completely positive matrices using [8], Proposition 2.2. Since the set of completely positive matrices is closed (cf. [8], Theorem 2.2), this implies that $\text{var}(L_1)$ is completely positive.

For part (b) let $k \in \mathbb{N}$ and $B \in M_{d,k}(\mathbb{R})$ with only non-negative entries be such that $C = BB^{\text{T}}$. Let further $(L_{i,t})_{t \in \mathbb{R}^+}$ with $i = 1, 2, \ldots, k$ be $k$ independent univariate subordinators with finite second moments and $\text{var}(\tilde{L}_1) = I_k$ where $\tilde{L}_t = (L_{1,t}, L_{2,t}, \ldots, L_{k,t})^{\text{T}}$ for $t \in \mathbb{R}^+$. Then $L := (B\tilde{L}_t)_{t \in \mathbb{R}^+}$ is a $d$-dimensional multivariate subordinator with $\text{var}(L_1) = BB^{\text{T}} = C$.

Combining (a) and (b) gives (c) and, finally, (d) follows from (c), because for $d \leq 4$ a matrix is doubly non-negative if and only if it is completely positive (cf. [21] or [8], Theorem 2.4). □

Part (d) provides a very nice complete characterisation of the covariance matrices of multivariate subordinators for $d \leq 4$. However, it does not extend to higher dimensions, since for dimensions five and greater there exist examples of doubly non-negative matrices that are not completely positive (see, e.g., [7]). In general dimensions it should thus be noted that there are easy-to-check, sufficient conditions for complete positivity, but the necessary and sufficient conditions known until now are more involved (cf. [36]).

To conclude this discussion of the second order properties of multivariate subordinators, we now strengthen part (b) of the last theorem by also considering the expected value.

**Proposition 5.3.** *Let $C \in M_d(\mathbb{R})$ be completely positive and $\mu \in \mathbb{R}^d$ have only strictly positive entries. Then there exists a $d$-dimensional multivariate subordinator $L$ such that $E(L_1) = \mu$ and $\text{var}(L_1) = C$.*

**Proof.** W.l.o.g. assume $C \neq 0$. (In the case $C = 0$, simply take $L$ as the deterministic Lévy process with drift $\mu$.) Let $k \in \mathbb{N}$ and $B \in M_{d,k}(\mathbb{R})$ with only non-negative entries be such that $C = BB^{\text{T}}$. Denote by $\mathbf{e} \in \mathbb{R}^k$ the vector $(1, 1, \ldots, 1)^{\text{T}}$ and set

$$\lambda = \min_{i=1,2,\ldots,d}\left\{\frac{\mu_i}{(B\mathbf{e})_i}\right\} \quad \text{and} \quad r = \frac{\lambda^3}{2},$$

where the well-definedness of the minimum follows from $C \neq 0$. Defining now $\gamma_L = \mu - 0.5\lambda B\mathbf{e}$ the choice of $\lambda$ ensures that $\gamma_L$ has only strictly positive entries.

Let further $(L_{i,t})_{t \in \mathbb{R}^+}$ with $i = 1, 2, \ldots, k$ be $k$ independent univariate compound Poisson processes with rate $r$ and the jump distribution being the exponential distribution with parameter $\lambda$. Define a $k$-dimensional Lévy process $\tilde{L}$ by $\tilde{L}_t = (L_{1,t}, L_{2,t}, \ldots, L_{k,t})^{\text{T}}$ for $t \in \mathbb{R}^+$. Then elementary calculations imply $E(\tilde{L}_1) = 0.5\lambda\mathbf{e}$ and $\text{var}(\tilde{L}_1) = I_k$.

Hence, defining $L$ by $L_t = \gamma_L t + B\tilde{L}_t$ for $t \in \mathbb{R}^+$ concludes the proof. □

For a general matrix subordinator $L$ we know that

$$\text{var}(\text{vec}(L_1)) = \int_{\mathbb{S}_d^+} \text{vec}(x)\text{vec}(x)^{\text{T}} \nu_L(\text{d}x),$$

but this does not seem to have any nice linear algebraic implications. Yet, the diagonal of any matrix subordinator forms a multivariate subordinator, so that the above results allow us at least to characterise the second order structure of the diagonal elements.

## 5.2. Normal mixture quadratic variation matrix subordinators

In this section we establish that for two special classes of quadratic variation matrix subordinators, which are special normal mixtures and hence in particular have elliptical distributions, one gets a particularly nice second order structure because only very well-understood matrices are involved.

Let us assume for a moment that $L$ is a $d$-dimensional compound Poisson process with rate one and the jump distribution being the $d$-dimensional normal distribution with mean zero and covariance matrix $C \in \mathbb{S}_d^+$. This implies that the discontinuous part of the quadratic variation $[L,L]^{\mathfrak{d}}$, which is given by

$$[L,L]^{\mathfrak{d}}_t = \sum_{0 < s \leq t} \Delta L_s (\Delta L_s)^{\mathrm{T}} \qquad \forall t \in \mathbb{R}^+,$$

is a compound Poisson process with rate one where the jump distribution is a Wishart distribution. Then, denoting the $d$-dimensional normal distribution with vanishing mean and variance $C$ by $N_C(\mathrm{d}x)$ and noting that $\operatorname{vec}(xx^{\mathrm{T}})\operatorname{vec}(xx^{\mathrm{T}})^{\mathrm{T}} = (x \otimes x)(x^{\mathrm{T}} \otimes x^{\mathrm{T}}) = (xx^{\mathrm{T}}) \otimes (xx^{\mathrm{T}})$ we have

$$E([L,L]^{\mathfrak{d}}_1) = \int_{\mathbb{R}^d} xx^{\mathrm{T}} \nu_L(\mathrm{d}x) = \int_{\mathbb{R}^d} xx^{\mathrm{T}} N_C(\mathrm{d}x) = C, \tag{5.1}$$

$$\operatorname{var}(\operatorname{vec}([L,L]^{\mathfrak{d}}_1)) = \int_{\mathbb{R}^d} \operatorname{vec}(xx^{\mathrm{T}})\operatorname{vec}(xx^{\mathrm{T}})^{\mathrm{T}} \nu_L(\mathrm{d}x) = \int_{\mathbb{R}^d} (xx^{\mathrm{T}}) \otimes (xx^{\mathrm{T}}) N_C(\mathrm{d}x) \tag{5.2}$$

$$= (C^{1/2} \otimes C^{1/2}) \int_{\mathbb{R}^d} (xx^{\mathrm{T}}) \otimes (xx^{\mathrm{T}}) N_{I_d}(\mathrm{d}x) (C^{1/2} \otimes C^{1/2})$$

$$= (C^{1/2} \otimes C^{1/2})(I_{d^2} + K_d + \operatorname{vec}(I_d)\operatorname{vec}(I_d)^{\mathrm{T}})(C^{1/2} \otimes C^{1/2})$$

$$= C \otimes C + K_d(C \otimes C) + \operatorname{vec}(C)\operatorname{vec}(C)^{\mathrm{T}} \tag{5.3}$$

from [18], Theorems 4.1 and 3.2(xii). Here $K_d \in M_{d^2}(\mathbb{R})$ denotes the commutation matrix which can be characterised by $K_d \operatorname{vec}(A) = \operatorname{vec}(A^{\mathrm{T}})$ for all $A \in M_d(\mathbb{R})$ (see [18] for more details). This can be easily generalised to the following result:

**Proposition 5.4.** *Let $L$ be a $d$-dimensional compound Poisson process with rate $r$ and the jumps being distributed like $(\varepsilon C)^{1/2} X$ where $X$ is a $d$-dimensional standard normal random variable, $C \in \mathbb{S}_d^+$ and $\varepsilon$ is a random variable in $\mathbb{R}^+$ with finite variance and independent of $X$. Then*

$$E([L,L]^{\mathfrak{d}}_1) = rE(\varepsilon)C, \tag{5.4}$$

$$\operatorname{var}(\operatorname{vec}([L,L]^{\mathfrak{d}}_1)) = rE(\varepsilon^2)(C \otimes C + K_d(C \otimes C) + \operatorname{vec}(C)\operatorname{vec}(C)^{\mathrm{T}}). \tag{5.5}$$

Moving away from a Lévy process of finite activity, a similar result holds for the following variant of type $G$ processes. Again, they correspond to a special kind of a normal mixture law such that infinite divisibility is ensured.

**Definition 5.5 (Type $\overline{G}$).** *Let $L$ be a d-dimensional Lévy process. If there exists an $\mathbb{R}^+$-valued infinitely divisible random variable $\varepsilon$ independent of a d-dimensional standard normal random variable $X$ and a matrix $C \in \mathbb{S}_d^+$ such that $L_1 \stackrel{\mathscr{L}}{=} (\varepsilon C)^{1/2} X$, then $L$ is said to be of* type $\overline{G}$. *(Here $\stackrel{\mathscr{L}}{=}$ denotes equality in law.)*

We have chosen the term "type $\overline{G}$" above because these processes correspond to a particular case of multG laws as defined in [2], Definition 3.1. Actually, many interesting Lévy processes are of type $\overline{G}$; for instance, the multivariate non-skewed and centred generalised hyperbolic (or normal inverse Gaussian) ones (cf. [10, 26]). For details on distributions/Lévy processes of type $G$ in general, we refer to [2, 17] and the references therein.

**Proposition 5.6.** *Let $L$ be a d-dimensional Lévy process of type $\overline{G}$ with a finite fourth moment. Then*

$$E([L,L]_1^{\mathfrak{d}}) = E(\varepsilon)C, \tag{5.6}$$

$$\operatorname{var}(\operatorname{vec}([L,L]_1^{\mathfrak{d}})) = \operatorname{var}(\varepsilon)(C \otimes C + K_d(C \otimes C) + \operatorname{vec}(C)\operatorname{vec}(C)^{\mathrm{T}}). \tag{5.7}$$

**Proof.** Let $\varepsilon$ be as in the definition of type $\overline{G}$ and $\nu_\varepsilon$ be its Lévy measure. Then by [2], Remark 3.1, $L$ has Lévy measure $\nu_L(\mathrm{d}x) = \int_{\mathbb{R}^+} N_{\tau C}(\mathrm{d}x)\nu_\varepsilon(\mathrm{d}\tau)$ and the finiteness of the fourth moment of $L$ implies that $E(\varepsilon^2) < \infty$ (unless $C = 0$, which is a trivial case). Hence,

$$E([L,L]_1^{\mathfrak{d}}) = \int_{\mathbb{R}^d} xx^{\mathrm{T}}\nu_L(\mathrm{d}x) = \int_{\mathbb{R}^+}\int_{\mathbb{R}^d} xx^{\mathrm{T}} N_{\tau C}(\mathrm{d}x)\nu_\varepsilon(\mathrm{d}\tau)$$

$$= C\int_{\mathbb{R}^+} \tau\nu_\varepsilon(\mathrm{d}\tau) = E(\varepsilon)C,$$

$$\operatorname{var}(\operatorname{vec}([L,L]_1^{\mathfrak{d}})) = \int_{\mathbb{R}^d}(xx^{\mathrm{T}}) \otimes (xx^{\mathrm{T}})\nu_L(\mathrm{d}x) = \int_{\mathbb{R}^+}\int_{\mathbb{R}^d}(xx^{\mathrm{T}}) \otimes (xx^{\mathrm{T}}) N_{\tau C}(\mathrm{d}x)\nu_\varepsilon(\mathrm{d}\tau)$$

$$= \int_{\mathbb{R}^+} \tau^2 \nu_\varepsilon(\mathrm{d}\tau)(C \otimes C + K_d(C \otimes C) + \operatorname{vec}(C)\operatorname{vec}(C)^{\mathrm{T}})$$

$$= \operatorname{var}(\varepsilon)(C \otimes C + K_d(C \otimes C) + \operatorname{vec}(C)\operatorname{vec}(C)^{\mathrm{T}})$$

using [18], Theorem 4.3, and noting that $\int_{\mathbb{R}^+} \tau\nu_\varepsilon(\mathrm{d}\tau)$ exists, because $\varepsilon$ is infinitely divisible and concentrated on $\mathbb{R}^+$ (cf. [30]). □

**Remark 5.7.** In the set-up of Proposition 5.6 $C$ and $E(\varepsilon), \operatorname{var}(\varepsilon)$ cannot be identified from the first and second moment of $[L,L]^{\mathfrak{d}}$. However, they are identified as soon as one fixes the scale of either $C$ or $\varepsilon$, for example, by demanding $E(\varepsilon) = 1$, $\rho(C) = 1, \|C\| = 1$ or $\det(C) = 1$.

Likewise, in the set-up of Proposition 5.4 one cannot identify $C$, $E(\varepsilon), E(\varepsilon^2)$ and the rate $r$ unless one fixes the scale of two of $C, \varepsilon$ and $r$.

Finally, let us consider a concrete example, viz. when $L$ is a $d$-dimensional centred and non-skewed generalised hyperbolic (GH) or normal inversion Gaussian Lévy (NIG) process (see [5, 10, 20, 26] for a comprehensive list of references regarding the univariate case). Note that the NIG Lévy processes form a subclass of the GH ones by fixing the parameter $\nu = -1/2$. The processes belonging to this subclass are highly tractable. Let $L$ be a $\mathrm{GH}(\nu, \alpha, 0, 0, \delta, C)$ Lévy process with $\nu \in \mathbb{R}$, $\alpha, \delta \in \mathbb{R}^{++}$ and $C \in \mathbb{S}_d^+$, that is, the location parameter $\mu$ and the skewness parameter $\beta$ are fixed to zero. Furthermore, let $\varepsilon$ be a generalised inverse Gaussian $\mathrm{GIG}(\nu, \delta, \alpha)$ random variable (see [14]) independent of a $d$-dimensional standard normal variable $X$, then

$$L_1 \stackrel{\mathscr{L}}{=} (\varepsilon C)^{1/2} X$$

and so $L$ is clearly type $\overline{G}$. From [14], page 13, or [5], Lemma 4, we obtain

$$E(\varepsilon) = \frac{\delta}{\alpha} \frac{K_{\nu+1}(\delta\alpha)}{K_\nu(\delta\alpha)}, \qquad \mathrm{var}(\varepsilon) = \frac{\delta^2}{\alpha^2} \frac{K_{\nu+2}(\delta\alpha) K_\nu(\delta\alpha) - K_{\nu+1}(\delta\alpha)^2}{K_\nu(\delta\alpha)^2},$$

where $K_\nu$ is the modified Bessel function of the third kind and order $\nu$ (see, e.g., [12, 35]), which can be defined via $K_\nu(z) = \frac{1}{2}\int_0^\infty y^{\nu-1} \mathrm{e}^{-z(y+y^{-1})/2}\,\mathrm{d}y$ for $z \in \mathbb{R}^{++}$. Hence, we obtain completely explicit representations for $E([L, L]_1^{\mathrm{o}})$ and $\mathrm{var}(\mathrm{vec}([L, L]_1^{\mathrm{o}}))$ in (5.6) and (5.7). If $L$ is actually an NIG Lévy process the representations simplify further, because in this case

$$E(\varepsilon) = \frac{\delta}{\alpha}, \qquad \mathrm{var}(\varepsilon) = \frac{\delta}{\alpha^3}$$

using that $K_\nu(z) = K_{-\nu}(z)$ and $K_{3/2}(z) = K_{1/2}(z)(1 + z^{-1})$.

Usually one assumes $\det(C) = 1$ in applications to make the parameters of the GH (NIG) distribution identified. We have opted not to make this assumption, because the above results are also of interest for degenerate GH Lévy processes in which $C$ is not invertible.

## Acknowledgement

Helpful comments of the anonymous referee are gratefully acknowledged.

## References

[1] Barndorff-Nielsen, O.E., Jensen, J.L. and Sørensen, M. (1998). Some stationary processes in discrete and continuous time. *Adv. in Appl. Probab.* **30** 989–1007. MR1671092

[2] Barndorff-Nielsen, O.E. and Pérez-Abreu, V. (2002). Extensions of type G and marginal infinite divisibility. *Theory Probab. Appl.* **47** 202–218. MR2001835

[3] Barndorff-Nielsen, O.E. and Pérez-Abreu, V. (2008). Matrix subordinators and related Upsilon transformations. *Theory Probab. Appl.* **52** 1–23. MR2354571

[4] Barndorff-Nielsen, O.E. and Shephard, N. (2001). Non-Gaussian Ornstein-Uhlenbeck-based models and some of their uses in financial economics (with discussion). *J. R. Stat. Soc. Ser. B Stat. Methodol.* **63** 167–241. MR1841412

[5] Barndorff-Nielsen, O.E. and Stelzer, R. (2005). Absolute moments of generalized hyperbolic distributions and approximate scaling of normal inverse Gaussian Lévy processes. *Scand. J. Statist.* **32** 617–637. MR2232346

[6] Barndorff-Nielsen, O.E. and Stelzer, R. (2007). Positive-definite matrix processes of finite variation. *Probab. Math. Statist.* **27** 3–43. MR2353270

[7] Berman, A. (1988). Complete positivity. *Linear Algebra Appl.* **107** 57–63. MR0960137

[8] Berman, A. and Shaked-Monderer, N. (2003). *Completely Positive Matrices*. Singapore: World Scientific. MR1986666

[9] Bhatia, R. (1997). *Matrix Analysis*, *Graduate Texts in Mathematics* **169**. New York: Springer. MR1477662

[10] Blæsild, P. and Jensen, J.L. (1981). Multivariate distributions of hyperbolic type. In *Statistical Distributions in Scientific Work – Proceedings of the NATO Advanced Study Institute Held at the Université degli Studi di Trieste, Triest, Italy, July 10–August 1, 1980* (C. Taillie, G.P. Patil and B.A. Baldessari, eds.) **4** 45–66. Dordrecht: Reidel. MR0656144

[11] Choi, M.-D. (1975). Positive semidefinite biquadratic forms. *Linear Algebra Appl.* **12** 95–100. MR0379365

[12] Gradshteyn, I. and Ryzhik, I.W. (1965). *Tables of Integrals, Series and Products*. New York: Academic Press. MR0197789

[13] Horn, R.A. and Johnson, C.R. (1991). *Topics in Matrix Analysis*. Cambridge: Cambridge Univ. Press. MR1091716

[14] Jørgensen, B. (1982). *Statistical Properties of the Generalized Inverse Gaussian Distribution*. New York: Springer. MR0648107

[15] Jurek, Z.J. and Mason, D.J. (1993). *Operator-limit Distributions in Probability Theory*. New York: Wiley. MR1243181

[16] Li, C.-K. and Pierce, S. (2001). Linear preserver problems. *Amer. Math. Monthly* **108** 591–605. MR1862098

[17] Maejima, M. and Rosiński, J. (2002). Type $G$ distributions on $\mathbb{R}^d$. *J. Theoret. Probab.* **15** 323–341. MR1898812

[18] Magnus, J.R. and Neudecker, H. (1979). The commutation matrix: Some properties and applications. *Ann. Statist.* **7** 381–394. MR0520247

[19] Marquardt, T. and Stelzer, R. (2007). Multivariate CARMA processes. *Stochastic Process. Appl.* **117** 96–120. MR2287105

[20] Masuda, H. (2004). On multidimensional Ornstein–Uhlenbeck processes driven by a general Lévy process. *Bernoulli* **10** 97–120. MR2044595

[21] Maxfield, J.E. and Minc, H. (1962). On the matrix equation $X'X = A$. *Proc. Edinburgh Math. Soc.* **13** 125–129. MR0174570

[22] Métivier, M. (1982). *Semimartingales: A Course on Stochastic Processes*, *De Gruyter Studies in Mathematics* **2**. Berlin: Walter de Gruyter. MR0688144

[23] Métivier, M. and Pellaumail, J. (1980). *Stochastic Integration*. New York: Academic Press. MR0578177

[24] Pierce, S., Lim, M.H., Loewy, R., Li, C.K., Tsing, N.K., McDonald, B. and Beasley, L. (1992). A survey of linear preserver problems. *Linear Multilinear Algebra* **33** 1–129. MR1346777

[25] Pigorsch, C. and Stelzer, R. (2009). A multivariate Ornstein–Uhlenbeck type stochastic volatility model. Available at http://www-m4.ma.tum.de.

[26] Prause, K. (1999). The generalized hyperbolic model: Estimation, financial derivatives and risk measures. Ph.D. thesis, Mathematische Fakultät, Albert–Ludwigs-Universität Freiburg i. Br., Freiburg, Germany.

[27] Protter, P. (2004). *Stochastic Integration and Differential Equations*, 2nd ed. *Stochastic Modelling and Applied Probability* **21**. New York: Springer. MR2020294

[28] Rajput, B.S. and Rosinski, J. (1989). Spectral representations of infinitely divisible processes. *Probab. Theory Related Fields* **82** 451–487. MR1001524

[29] Rudin, W. (1976). *Principles of Mathematical Analysis*, 3rd ed. Singapore: McGraw-Hill. MR0385023

[30] Sato, K. (1999). *Lévy Processes and Infinitely Divisible Distributions. Cambridge Studies in Advanced Mathematics* **68**. Cambridge Univ. Press. MR1739520

[31] Sato, K. (2006). Addtitive processes and stochastic integrals. *Illinois J. Math.* **50** 825–851. MR2247848

[32] Sato, K. and Yamazato, M. (1984). Operator-selfdecomposable distributions as limit distributions of processes of Ornstein–Uhlenbeck type. *Stochastic Process. Appl.* **17** 73–100. MR0738769

[33] Schneider, H. (1965). Positive operators and an inertia theorem. *Numer. Math.* **7** 11–17. MR0173678

[34] Skorohod, A.V. (1991). *Random Processes with Independent Increments. Mathematics and Its Applications (Soviet Series)* **47**. Dordrecht: Kluwer. MR1155400

[35] Watson, G.N. (1952). *A Treatise on the Theory of Bessel Functions*, reprinted 2nd ed. Cambridge Univ. Press. MR0010746

[36] Xu, C. (2004). Completely positive matrices. *Linear Algebra Appl.* **379** 319–327. MR2039746

[37] Yamazato, M. (1983). Absolute continuity of operator-selfdecomposable distributions on $\mathbb{R}^d$. *J. Multivariate Anal.* **13** 550–560. MR0727040